\documentclass[letterpaper,11pt]{amsart}

\usepackage[margin=2cm]{geometry}
\usepackage{amscd, amssymb}
\usepackage{amsmath,amscd}
\usepackage[driverfallback=hypertex]{hyperref}
\usepackage{graphics}
\usepackage{graphicx}
\usepackage{color}
\usepackage{bm}
\usepackage{enumitem}
\setlist[enumerate,1]{label={(\alph*)}}
\setlist[enumerate,2]{label={(\roman*)}}
\input xy
\xyoption{all}

\newtheorem{thm}{Theorem}[section]

\newcommand{\ignore}[1]{}

\usepackage{color}

\definecolor{Red}{rgb}{1,0,0}
\definecolor{Blue}{rgb}{0,0,1}
\definecolor{Olive}{rgb}{0.41,0.55,0.13}
\definecolor{Yarok}{rgb}{0,0.5,0}
\definecolor{Green}{rgb}{0,1,0}
\definecolor{MGreen}{rgb}{0,0.8,0}
\definecolor{DGreen}{rgb}{0,0.55,0}
\definecolor{Yellow}{rgb}{1,1,0}
\definecolor{Cyan}{rgb}{0,1,1}
\definecolor{Magenta}{rgb}{1,0,1}
\definecolor{Orange}{rgb}{1,.5,0}
\definecolor{Violet}{rgb}{.5,0,.5}
\definecolor{Purple}{rgb}{.75,0,.25}
\definecolor{Brown}{rgb}{.75,.5,.25}
\definecolor{Grey}{rgb}{.5,.5,.5}

\begin{document}
\title{A Cayley-type identity for trees}
\author{Ran J. Tessler}
\address{Institute for Theoretical Studies, ETH Z\"urich}

\begin{abstract}
We prove a weighted generalization of the formula for the number of plane vertex-labeled trees.
\end{abstract}

\maketitle

\section{Introduction}
It is well known that the number of vertex-labeled trees on $n$ vertices is $n^{n-2}.$ The formula was discovered by Carl Wilhelm Borchardt in 1860 \cite{Borch} and was extended by Cayley in \cite{Cayley2}. Since then many proofs of this formula were given, and many extensions were found. A beautiful well-known extension is the following weighted Cayley formula.
\begin{thm}
Let $T_n$ be the set of vertex-labeled trees with $n$ vertices labeled by $[n]=\{1,\ldots,n\}$.
Associate a variable $x_i$ to every $i \in [n]$, and associate the monomial $\prod_{i \in [n]}x_i^{d_T(i)}$ to $T\in T_n$, where $d_T(i)$ is the degree of $i$ in $T.$ Then
\[
\sum_{T\in T_n} \prod_{i \in [n]}x_i^{d_T(i)} = \prod_{i=1}^n x_i\left(\sum_{i =1}^n{x_i}\right)^{n-2}.
\]
\end{thm}
An identity closely related to Cayley's formula is the formula, due to Leroux and Miloudi, \cite{LerMil} (see also \cite{Callan} for a short proof) which says that for $n\geq2$ there are $\binom{2n-3}{n-1}$ vertex-labeled plane trees on $n$ vertices. By a plane tree we mean an abstract tree enriched with cyclic orders for the edges which emanate from each vertex.

In this note we prove a "weighted version" for this formula, namely
\begin{thm}
Associate a variable $x_i$ to every $i \in [n],$ and for integers $m\geq 1$ and $a$ denote by $\binom{x+m+a}{x+a}$ the polynomial \[\prod_{i=a+1}^{a+m}(x+i),\]we extend the definition to $m=0$ by writing $\binom{x+a}{x+a}=1.$ For $n\geq 2$ it holds that
\[\sum_{T\in T_n}\prod_{i \in [n]}\binom{x_i+d_T(i)-1}{x_i-1} = \prod_{i=1}^n x_i\binom{\sum_{i \in [n]}{x_i}+2n-3}{\sum_{i \in [n]}{x_i} + n-1}.
\]
\end{thm}
For example, for $n=2$ the left and right hand sides of the formula give $x_1x_2.$ For $n=3$ the formula gives \[x_1x_2x_3(x_1+x_2+x_3+3).\]
Dividing both sides by $\prod x_i$ and substituting $x_1=\ldots=x_n=0$ gives precisely the Leroux-Miloudi formula.
The weighted Leroux-Miloudi formula was used in \cite{Bomba} to calculate and prove the threshold for the appearance of spanning $2-$spheres in the Linial-Meshulam model for random $2-$complexes.
\subsection{Acknowledgements}
R.T. is supported by Dr. Max R\"ossler, the Walter Haefner Foundation and the ETH Z\"urich
Foundation.
\section{Proof of the formula}
When $n=2$ the formula trivially holds (in fact, correctly interpreted, the formula extends to $n=1$). Our proof will be inductive.
By dividing both sides by $\prod x_i$ the theorem is seen to be equivalent to proving, for $n\geq 2,$
\begin{equation}\label{eq:1}
\sum_{T\in T_n}\prod_{i \in [n]}\binom{x_i+d_T(i)-1}{x_i} =\binom{\sum_{i \in [n]}{x_i}+2n-3}{\sum_{i \in [n]}{x_i} + n-1}.
\end{equation}
Substitute $y_i=x_i+1.$ \eqref{eq:1} then translates to
\begin{equation}\label{eq:2}
L_n(y_1,\ldots,y_n):=\sum_{T\in T_n}\prod_{i \in [n]}\binom{y_i+d_T(i)-2}{y_i-1} =\binom{\sum_{i \in [n]}{y_i}+n-3}{\sum_{i \in [n]}{y_i} -1}=:R_n(y_1,\ldots,y_n).
\end{equation}
Both the left hand and right hand side are polynomials of degree $n-2$ in $n$ variables. Thus, any monomial does not contain at least one of the variables. Hence, \eqref{eq:2} will follow from proving that for each $i=1,\ldots, n$
\begin{equation}\label{eq:3}
L_n(y_1,\ldots,y_n)|_{y_i=0}=R_n(y_1,\ldots,y_n)|_{y_i=0}.
\end{equation}
Since $L_n,R_n$ are in addition symmetric, it is enough to prove \eqref{eq:3} for $i=n.$
As \[\left(\sum_{i=1}^ny_i\right)|_{y_n=0}=\sum_{i=1}^{n-1}y_i\] we have
\begin{equation}\label{eq:4}
R_n|_{y_n=0}=(n-3+\sum_{i=1}^{n-1}y_i)R_{n-1}.
\end{equation}
The induction will therefore follow if we could show that
\begin{equation}\label{eq:5}
L_n|_{y_n=0}=(n-3+\sum_{i=1}^{n-1}y_i)L_{n-1}.
\end{equation}
Denote by $w(T)$ the summand in \eqref{eq:2} which corresponds to the tree $T,$
\[w(T)=\prod_{i \in [n]}\binom{y_i+d_T(i)-2}{y_i-1}.\]
Observe that if $y_n=0$ then $w(T)=0$ whenever $d_T(n)>1.$
Thus,
\begin{equation}\label{eq:6}
L_n(y_1,\ldots,y_n)|_{y_i=0}=\sum_{T\in T'_n}\prod_{i \in [n]}\binom{y_i+d_T(i)-2}{y_i-1},
\end{equation}
where $T'_n\subseteq T_n$ is the collection of trees in which $n$ is a leaf.
For $t\in T'_n$ let $a(t)\in[n-1]$ be the single neighbour of $n$ and let $t(T)\in T_{n-1}$ be the tree obtained from erasing the vertex $n.$
It hold that
\begin{equation}\label{eq:7}
w(T)= (y_{a(t)}+d_T(a(t))-2)w(t(T))=(y_{a(t)}+d_{t(T)}(a(t))-1)w(t(T)).
\end{equation}
Since the sum of degrees of vertices in a graph is twice the number of edges, and the number of edges in a tree on $m$ vertices is $m-1$ \eqref{eq:7} yields, for any $T\in T_{n-1},$
\begin{equation}\label{eq:8}
\sum_{T'\in t^{-1}(T)}w(T')= \sum_{a\in[n-1]}(y_{a}+d_{T}(a)-1)w(T)=(\sum_{a=1}^{n-1}y_i+n-3)w(T).
\end{equation}
Putting \eqref{eq:6},\eqref{eq:8} and the definition of $L_n$ together
\[L_n|_{y_i=0}=\sum_{T\in T'_n}\prod_{i \in [n]}\binom{y_i+d_T(i)-2}{y_i-1}=(\sum_{a=1}^{n-1}y_i+n-3)\sum_{T\in T_{n-1}}w(T)=L_{n-1}\]
which is precisely \eqref{eq:5}
\bibliographystyle{abbrvnat}
\bibliography{bibli_comb}

\end{document}